\documentclass[a4paper]{article}
\usepackage{fullpage}
\usepackage{amsmath}
\usepackage{amsthm}
\usepackage{amssymb}
\usepackage{amsbsy}
\usepackage{mathrsfs}
\usepackage{amstext,hyperref}  

\usepackage{chngcntr}
\usepackage{changepage}

\usepackage{array,longtable}
\usepackage{enumerate}

\usepackage{fullpage} 

\newtheorem{theorem}{Theorem}[section]

\newtheorem{lemma}[theorem]{Lemma} 
\newtheorem{corollary}[theorem]{Corollary}

\newtheorem{conjecture}{Conjecture}[theorem] 
\newtheorem{observation}{Observation}[theorem] 
\newtheorem{remark}{Remark}[theorem] 

\newtheorem*{thm-main}{Theorem~\ref{thm-main}}

\hypersetup{
	colorlinks=true, 
	linkcolor=red,
	citecolor=blue,
}

\usepackage{graphicx,accents}

\makeatletter
\DeclareRobustCommand{\ovl}[1]{%
	\mathpalette\do@cev{#1}%
}
\newcommand{\do@cev}[2]{%
	\fix@cev{#1}{+}%
	\reflectbox{$\m@th#1\ovr{\reflectbox{$\fix@cev{#1}{-}\m@th#1#2\fix@cev{#1}{+}$}}$}%
	\fix@cev{#1}{-}%
}
\newcommand{\fix@cev}[2]{%
	\ifx#1\displaystyle
	\mkern#23mu
	\else
	\ifx#1\textstyle
	\mkern#23mu
	\else
	\ifx#1\scriptstyle
	\mkern#22mu
	\else
	\mkern#22mu
	\fi
	\fi
	\fi
}
\usepackage[title]{appendix}
\usepackage{soul}
\usepackage{mathtools}  
\usepackage{pgfplots}
\usepackage{tikz}
\usepackage{ifthen}
\usetikzlibrary{shapes.geometric, positioning, calc}
\tikzset{%
	mynode/.style={%
		circle, minimum size=1.21mm, inner sep=0pt, draw=black, fill=black
	}
}

\usepackage{todonotes}

\definecolor{lightblue}{rgb}{0.68,0.85,0.9}

\begin{document}
	
	\title{{\bf Large planar $(n,m)$-cliques}}
	
	\author{{\sc Susobhan Bandopadhyay}$\,^{a}$, {\sc Sagnik Sen}$\,^{b}$, 
		 {\sc 
			S Taruni}{$\,^{c}$}
		\mbox{}\\
		{\small $(a)$ School of Technology and Computer Science,} \\ {\small Tata Institute of Fundamental Research, Mumbai India.}\\
		{\small $(b)$ Indian Institute of Technology Dharwad, India. }\\
		{\small $(c)$ Centro de Modelamiento Matemático (CNRS IRL2807),}\\
		{\small Universidad de Chile, Santiago, Chile.}\\
	}
	
	\date{}
	
	\maketitle

\begin{abstract}
An \textit{$(n,m)$-graph} $G$ is a graph 
having both arcs and edges, and its arcs (resp., edges) are labeled using one of the $n$ (resp., $m$)  different symbols. 
An \textit{$(n,m)$-complete graph} $G$ is an $(n,m)$-graph without loops or multiple edges in its underlying graph such that  
identifying any pair of vertices results in a loop or parallel adjacencies with distinct labels. 
We show that a planar $(n,m)$-complete graph cannot have more than $3(2n+m)^2+(2n+m)+1$ vertices, for all $(n,m) \neq (0,1)$ and that the bound is tight. 
This positively settles a conjecture by Bensmail \textit{et al.}~[Graphs and Combinatorics 2017]. 
 \end{abstract}

\section{Introduction and the main result}\label{sec intro}
All graphs considered in this article are finite and contains no loop or multiple edges. 

\medskip

\noindent \textbf{Definitions, notation, and terminology:}  An \textit{$(n,m)$-graph} $G$ is a graph having $n$ different types of arcs and $m$ different types of edges.
Given an $(n,m)$-graph $G$, each of its arcs is labeled by one of the $n$ (even) numbers from 
$\{2, 4, \ldots, 2n\}$ and each of its edges is labeled by one of the numbers from 
$\{2n+1, 2n+2, \ldots, 2n+m\}$. Given an arc $uv$ of type $2i$, we will call $vu$ as a \textit{reverse arc} of type $2i-1$. Also, 
in such a scenario, $v$ is a $2i$-neighbor of $u$ and, equivalently, $u$ is a $(2i-1)$-neighbor of $v$. 
Moreover, for an edge $uv$ of type $2n+j$, the edge $vu$ is also 
an edge of type $2n+j$. So,
$u$ and $v$ are $(2n+j)$-neighbors of each other. This means, each natural number from 
$A_{n,m} = \{1, 2, \ldots, 2n+m\}$ denotes a particular adjacency (arc, reverse arc, or edge) and its type. Given 
$\alpha \in \{1, 2, \ldots, 2n+m\}$ and a vertex $v$ of an $(n,m)$-graph $G$, the 
symbols $N^{\alpha}(v)$ and $d^{\alpha}(v)$ denote the set and number of all $\alpha$-neighbors of $v$ in $G$. We will reserve the lower case 
Greek alphabets $\alpha, \beta, \gamma$ etc. to denote adjacency and its type in $(n,m)$-graphs without declaring their domain (sometimes) henceforth.

A \textit{homomorphism} of an $(n,m)$-graph $G$ to another $(n,m)$-graph $H$ is a vertex mapping 
$f: V(G) \to V(H)$ such that 
if $uv$ is an arc (resp., reverse arc, edge) of type $\alpha$ in $G$, then $f(u)f(v)$ is also an 
arc (resp., reverse arc, edge) of type $\alpha$ in $H$. If there exists a homomorphism of $G$ to $H$, then we denote it by $G \rightarrow H$. 
An \textit{$(n,m)$-complete graph} $G$ is an $(n,m)$-graph
which does not admit any homomorphism to an $(n,m)$-graph having less number of vertices. 
In other words, an $(n,m)$-complete graph $G$ is an $(n,m)$-graph in which identifying any pair of adjacent vertices results in a loop or multiple edges having distinct adjacency types.  The parameters, closely associated with homomorphisms, are listed below. 

Let $G$ be an $(n,m)$-graph. Then its $(n,m)$-chromatic number, denoted by 
$\chi_{n,m}(G)$ is given as 
$$\chi_{n,m}(G) = \min \{|V(H)|: G \rightarrow H\}.$$ 
Notice that, the definition can be equivalently reformulated with the restriction on $H$ to be an $(n,m)$-complete graph. 

An \textit{absolute $(n,m)$-clique} of $G$ is a vertex subset $A$ that induces an $(n,m)$-complete graph. 
The \textit{absolute $(n,m)$-clique number} of $G$, denoted by $\omega_{a(n,m)}(G)$, is given as 
$$\omega_{a(n,m)}(G) = \max \{|A|: A \text{ is an absolute $(n,m)$-clique of } G\}.$$

For a family $\mathcal{F}$ of graphs, the above two parameters can be defined as 
$$p(\mathcal{F}) = \max\{p(G) : \text{ the underlying graph of $G$ belongs to } \mathcal{F}\},$$
where 
$p \in \{\chi_{n,m}, \omega_{a(n,m)}\}$. 

The natural relation between the two parameters is clear from their definitions~\cite{bensmail2017analogues}: 
$$\omega_{a(n,m)}(\mathcal{F}) \leq \chi_{n,m}(\mathcal{F}).$$

\medskip

\noindent \textbf{Context, motivation, and our contributions:} Let $\mathcal{P}_3$ denote the family of planar graphs.
Thus, $\chi_{0,1}(\mathcal{P}_3)=4$ is precisely the Four-Color Theorem. However, finding the exact values of 
$\chi_{n,m}(\mathcal{P}_3)$ are open problems for
all $(n,m) \neq (0,1)$. In 1994 Courcelle~\cite{courcelle1994monadic} showed that every planar 
$(1,0)$-graph $G$ satisfying
$d^1(v) \leq 3$ for all $v \in V(G)$
have $\chi_{1,0}(G) \leq 4^33^{63}$. This bound was quickly improved by Raspaud and Sopena~\cite{raspaud1994good} by proving 
$\chi_{1,0}(\mathcal{P}_3) \leq 80$ using the fact that the acyclic chromatic number of planar graphs is at most $5$ due to Borodin~\cite{borodin1979acyclic}. For $(0,m)$-graphs, this bound was generalized by 
Alon and Marshall~\cite{alon1998homomorphisms} by proving 
$\chi_{0,m}(\mathcal{P}_3) \leq 5m^4$ while they provided a
lower bound of the same by a cubic polynomial of $m$. 
Later Gu\'{s}piel and Gutowski~\cite{guspiel2017universal} asymptotically improved the upper bound by proving 
$\chi_{0,m}(\mathcal{P}_3) = \Theta(m^3)$ (their upper bound beats the previous one only for very large values of $m$). 
Ne\v{s}et\v{r}il and Raspaud~\cite{nevsetvril2000colored}
generalized the former bound by showing 
$\chi_{n,m}(\mathcal{P}_3) \leq 5(2n+m)^4$. 
A lower bound given by a cubic polynomial of 
$(2n+m)$ for $\chi_{n,m}(\mathcal{P}_3)$ was established~\cite{fabila2008lower} by generalizing the works of Alon and Marshall~\cite{alon1998homomorphisms}. However, any improvement in the 
existing bounds of $\chi_{n,m}(\mathcal{P}_3)$ seems like an extremely difficult problem.

Attempts to partially solve the problem of finding $\chi_{n,m}(\mathcal{P}_3)$ is made by studying, and establishing lower and upper bounds of $\chi_{n,m}(\mathcal{F})$, where $\mathcal{F}$ represents different subfamilies of planar graphs~\cite{BORODIN2004147}. 
Especially, focused studies on this front are made 
for $(n,m)=(1,0)$ and $(0,2)$~\cite{naserasr2015homomorphisms, DBLP:journals/corr/OchemPS14, sopena2016homomorphisms}. 
In practice, almost all 
the lower bounds  are achieved using $(n,m)$-complete graphs, 
or their structural extensions, or 
through studying their interactions. Thus, the knowledge of 
planar $(n,m)$-complete graphs is a fundamentally important 
concept in the study of finding $\chi_{n,m}(\mathcal{P}_3)$.

Although the concept of $(n,m)$-complete graphs was 
known and implicitly used in the above-mentioned works, 
Klostermeyer and MacGillivray~\cite{klostermeyer2004analogues} formally introduced the 
$(1,0)$-complete graphs (under the name of \textit{ocliques})
alongside the parameter $\omega_{a(1,0)}$ (as the notation 
$\omega_o$)
and showed that $\omega_{a(1,0)}(\mathcal{P}_3) \leq 36$.  They~\cite{klostermeyer2004analogues} also conjectured 
$\omega_{a(1,0)}(\mathcal{P}_3) = 15$ which was positively resolved by 
Nandy \textit{et al.}~\cite{nandy2016outerplanar}. 
Later, this notion was generalized by Bensmail \textit{et al.}~\cite{bensmail2017analogues}. They~\cite{bensmail2017analogues} proved that $3(2n+m)^2+(2n+m)+1 \leq \omega_{a(n,m)}(\mathcal{P}_3) \leq 9(2n+m)^2+2(2n+m)+2 $. An equivalent formulation is as follows. 
\begin{theorem}[\cite{bensmail2017analogues}]
 If $G$ is a planar $(n,m)$-complete graph for all $(n,m) \neq (0,1)$, then $G$ has at most $9(2n+m)^2+2(2n+m)+2$ vertices. Moreover, there exists planar (n,m)-complete graphs having $3(2n+m)^2+(2n+m)+1$ vertices, for all $(n,m) \neq (0,1)$.   
\end{theorem}
And subsequently, they~\cite{bensmail2017analogues} conjectured that the lower bound is tight.
\begin{conjecture}[\cite{bensmail2017analogues}]\label{main conjecture}
    If $G$ is a planar $(n,m)$-complete graph for all $(n,m) \neq (0,1)$, then $G$ has at most $3(2n+m)^2+(2n+m)+1$ vertices. Moreover, the bound is tight.
\end{conjecture}

Our main contribution is 
to positively settle Conjecture~\ref{main conjecture}. Instead of finding the exact value of 
$\chi_{n,m}(\mathcal{P}_3)$, we are finding the exact value of a lesser quantity 
$\omega_{a(n,m)}(\mathcal{P}_3)$. Thus, one can consider solving this sub-question as a step towards finding the analog of the Four-color Theorem for $(n,m)$-graphs.

\begin{theorem}\label{th main}
    If $G$ is a planar $(n,m)$-complete graph for $(n,m) \neq (0,1)$, then $G$ has at most $3(2n+m)^2+(2n+m)+1$ vertices. Moreover, the bound is tight.  
\end{theorem}

\begin{remark}
  As pointed out in~\cite{guspiel2017universal, nevsetvril2000colored}, the $(n,m)$-graph analogs of extremely simple undirected graph problems are often found to be notoriously difficult.   Our problem (the statement of Theorem~\ref{th main}) also falls in the category where its undirected analog (finding the largest complete planar graph) is a classroom exercise problem, and the proof for the case when $(n,m) = (1,0)$ is the main result of~\cite{nandy2016outerplanar} whose proof spreads over $17$ pages. However, we cannot generalize this proof, as the cases handled in~\cite{nandy2016outerplanar} if applied to our proof, generates an exponential order of $(2n+m)$ many subcases.

\end{remark}

\noindent \textbf{Applications:}
Ne\v{s}et\v{r}il and Raspaud~\cite{nevsetvril2000colored} introduced the notion of ``homomorphisms of $(n,m)$-graphs'', under the name ``colored homomorphisms of colored mixed graphs'', as a generalization of homomorphisms of undirected graphs, oriented graphs, signed graphs, and $m$-edge-colored graphs~\cite{alon1998homomorphisms, hell2004graphs, sopena2016homomorphisms}. 
This generalization has connections with 
harmonious coloring of 
graphs~\cite{alon1998homomorphisms}, nowhere-zero 
flows~\cite{BORODIN2004147},
binary predicate logic~\cite{nevsetvril2000colored}, 
Coxeter groups~\cite{alon1998homomorphisms}. Moreover, 
from the application point of view, 
homomorphisms of $(n,m)$-graphs (and its variants) are natural model for the Query Evaluation Problem in graph database~\cite{angles2017foundations, 
beaudou2019complexity} 
which is widely used in
social networks (e.g., Facebook, Twitter), information networks (e.g., World Wide Web, citation of academic papers), technological networks (e.g., internet, Geographic Information Systems or \textit{GIS},
phone networks), and biological networks (e.g., genomics, food web, neural networks)~\cite{angles2008survey}. Furthermore, from theoretical point of view,  
$(n,m)$-complete graphs are studied from both structural and algorithmic angle, 
including for specific values of $(n,m)$~\cite{bensmail2017analogues, chakraborty2021clique,coelho2023absolute, dybizbanski2020oriented, klostermeyer2004analogues, nandy2016outerplanar,  naserasr2015homomorphisms}. 
Moreover, in the study of deeply critical oriented graphs~\cite{borodin2001deeply, DBLP:journals/jgt/DuffyDBS23}, all the examples found
to date are $(1,0)$-complete graphs in particular.

\medskip

\noindent \textit{Note:} For standard graph theoretic notation and terminology, please refer to West~\cite{west2001introduction}. Also, any standard term (degree, diameter, connected, planarity, etc.,) for undirected graph used in the context of  $(n,m)$-graph $G$ should be considered for its underlying graph $G$ instead.

\section{Proof of Theorem~\ref{th main}}
A \textit{special $2$-path} is a path in an $(n,m)$-graph $G$ of the form $uwv$ such that 
$w \in N^{\alpha}(u) \cap N^{\beta}(v)$
for some $\alpha \neq \beta$. 
A vertex \textit{$u$ sees $v$} if either they are adjacent, or they are connected by a 
special $2$-path. Moreover, if they are connected by a special $2$-path of the form $uwv$, then we say that \textit{$u$ sees $v$ through $w$}. 
For the proof, we will use the following handy characterization of 
$(n,m)$-complete graphs using the notion of special $2$-paths. 

\begin{observation}[\cite{bensmail2017analogues}]\label{obs char of (n,m)-complete}
    An $(n,m)$-graph $G$ is an $(n,m)$-complete graph if and only if any pair of vertices of $G$ sees each other. 
\end{observation}

Bensmail \textit{et al.}~\cite{bensmail2017analogues} provided examples of planar $(n,m)$-complete graphs on $3(2n+m)^2+(2n+m)+1$ vertices for all $(n,m)\neq (0,1)$. Thus, the ``moreover'' part of Theorem~\ref{th main} is proved. 
Additionally, Nandy, Sen, and Sopena~\cite{nandy2016outerplanar} 
proved Theorem~\ref{th main} for $(n,m) = (1,0)$, and the case $(n,m)=(0,2)$ can be handled similarly (a detailed proof can be found in~\cite{sen2014contribution}). 
Thus, given any planar $(n,m)$-complete graph $H$, it is enough to show that $H$ has at most $3(2n+m)^2+(2n+m)+1$ vertices, where the values of $(n,m)$ satisfies 
$2n+m \geq 3$. That means if $H$ has at most $3(3)^2+(3)+1 = 31$ vertices, then there is nothing to prove. However, using some known results, and some carefully chosen assumptions, it is enough to consider an $H$
having specific structural properties.

\subsubsection*{Initialization of $H$}  
First of all, notice that an $(n,m)$-complete graph remains an $(n,m)$-complete graph even if we add new arcs or edges (between non-adjacent vertices). Therefore, we may assume that $H$ is a triangulated planar graph. Moreover, we know that an $(n,m)$-complete graph has diameter at most $2$ due to Observation~\ref{obs char of (n,m)-complete}. 
As the largest planar graph having diameter $1$ is $K_4$, that is, the complete graph on $4$ vertices,  we can assume that the diameter of $H$ is equal to $2$.

It is known~\cite{goddard2002domination} that a planar graph having diameter $2$ has domination number at most $2$, except for one example on $11$ vertices. 
 As this exception trivially respects the upper bound,  we may assume that $H$ has domination number at most $2$. Bensmail \textit{et al.}~\cite{bensmail2017analogues} showed that if 
a planar $(n,m)$-complete graph has domination number
$1$, then it has at most $3(2n+m)^2+(2n+m)+1$ vertices. Therefore, we may further refine our assumption and suppose that $H$ has domination number exactly equal to $2$. 

Let $D = \{x,y\}$ be a dominating set of size $2$ in 
$H$ satisfying the following property: given any 
dominating set $D' = \{x',y'\}$ of size $2$ in $H$, 
we must have 
$$|N(x') \cap N(y')| \leq |N(x) \cap 
N(y)|.$$ 
That means, among all dominating sets of 
size $2$ in $H$, the vertices of $D$ has the maximum 
number of common neighbors. Let us fix this $D$ with the above-mentioned property for the rest of the proof. 

Therefore, $H$ is a triangulated plane 
$(n,m)$-complete graph with domination number $2$, and with a dominating set $D = \{x,y\}$ where 
$|N(x) \cap N(y)|$  is maximized. 

\subsubsection*{Set up of the proof}  
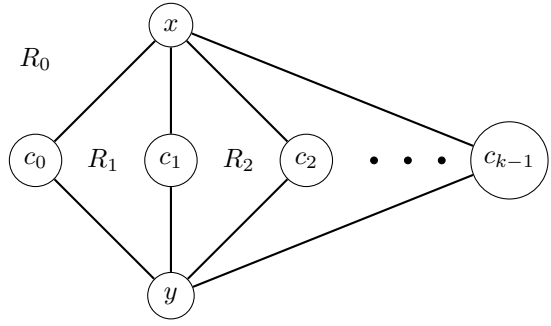
\begin{figure}
\begin{center}
		\begin{tikzpicture} 
        
			[scale=.7,auto=center] 

			\node[circle,draw] (a1) at (0,4) {$x$};
			\node[circle,draw] (a2) at (0,0) {$y$};
			\node[circle,draw] (a3) at (-2,2) {$c_0$};
			\node[circle,draw]  (a4) at (0,2) {$c_1$};
            \node[circle,draw] (a5) at (2,2) {$c_2$};
			\node[circle,draw]  (a6) at (5,2) {$c_{k-1}$};
			\node at (-2,3.5) {$R_0$};
			\node at (-1,2) {$R_1$};
            \node at (1,2) {$R_2$};
			\draw[fill] (3,2) circle[radius=0.05];
			\draw[fill] (3.5,2) circle[radius=0.05];
			\draw[fill] (4,2) circle[radius=0.05];
			\draw[thick] (a1) -- (a3);
			\draw[thick] (a1) -- (a4);
			\draw[thick] (a1) -- (a5);
			\draw[thick] (a1) -- (a6);
			\draw[thick] (a2) -- (a3);
			\draw[thick] (a2) -- (a4);
            \draw[thick] (a2) -- (a5);
            \draw[thick] (a2) -- (a6);

		\end{tikzpicture}

	\caption{A partial embedding of $H$.}
	\label{fig Fk}
    \end{center}
\end{figure}

Suppose $C = N(x) \cap N(y) = \{c_0, c_1, \ldots, c_{k-1}\}$ be the set of common neighbors of $x$ and $y$. Let us assume a (partial) planar embedding of $H$ where the vertices 
$c_0, c_1, \ldots, c_{k-1}$ are arranged in an anti-clockwise order around $x$, where $k \geq 2$. Let us, for now, just think about a planar embedding 
of the $(n,m)$-graph $K_{2,k}$ induced by the arcs, reverse arcs, and edges between the sets $C$ and $D$. 
Without loss of generality, we may assume that the face enclosed by the cycle $xc_0yc_{k-1}x$ is the outerface of this $K_{2,k}$. In the fixed embedding of $H$, we call this face the (unbounded) region $R_0$. Also, the face of the $K_{2,k}$ enclosed by the cycle $xc_{i-1}yc_{i}x$ is the region $R_i$ in 
the fixed embedding of $H$, where $i \in \{1, 2, \ldots, k-1\}$. If the boundaries of two distinct regions share a vertex of $C$, then they are \textit{adjacent regions}. See Fig.~\ref{fig Fk} for a pictorial reference.

Let $R$ be a region and $uv$ be an edge of its boundary. Let $Z = \{z_1, z_2, \ldots, z_t\}$  be a set of vertices contained in $R$ that are adjacent to both $u$ and $v$. Now, modify the plane graph by deleting all the vertices that are not part of the boundary of $R$ and $Z$. In this scenario, note that, it is possible to draw a line between the midpoint of the edge $uv$ and exactly one of the vertices from $Z$ without crossing any other edge or vertex. This vertex of $Z$, say $z_1$, is called the 
\textit{$1$-nearest} vertex of $Z$ to $uv$. Inductively, we can define \textit{$i$-nearest} neighbor inductively: if $z_1, z_2, \ldots, z_i$ are the $j$-nearest vertices of $Z$ to $uv$ for $j \in \{1, 2, \ldots, i\}$, then the $1$-nearest vertex of 
$Z \setminus \{z_1, z_2, \ldots, z_i\}$ to $uv$ is the $(i+1)$-nearest neighbor of $Z$ to $uv$. Moreover, the $i$-nearest neighbor is also called the \textit{$(t-i+1)$-farthest} neighbor of $Z$ to $uv$.

We define $C^{\alpha, \beta} = N^{\alpha}(x) \cap N^{\beta}(y)$ as the set of vertices that are $\alpha$-neighbors of $x$ and $\beta$-neighbors of $y$. Additionally, we denote $C_x^{\alpha} = C \cap N^{\alpha}(x)$ as the set of $\alpha$-neighbors of $x$ that are also adjacent to $y$, and $C_y^{\beta} = C \cap N^{\beta}(y)$ as the set of $\beta$-neighbors of $y$ that are also adjacent to $x$.

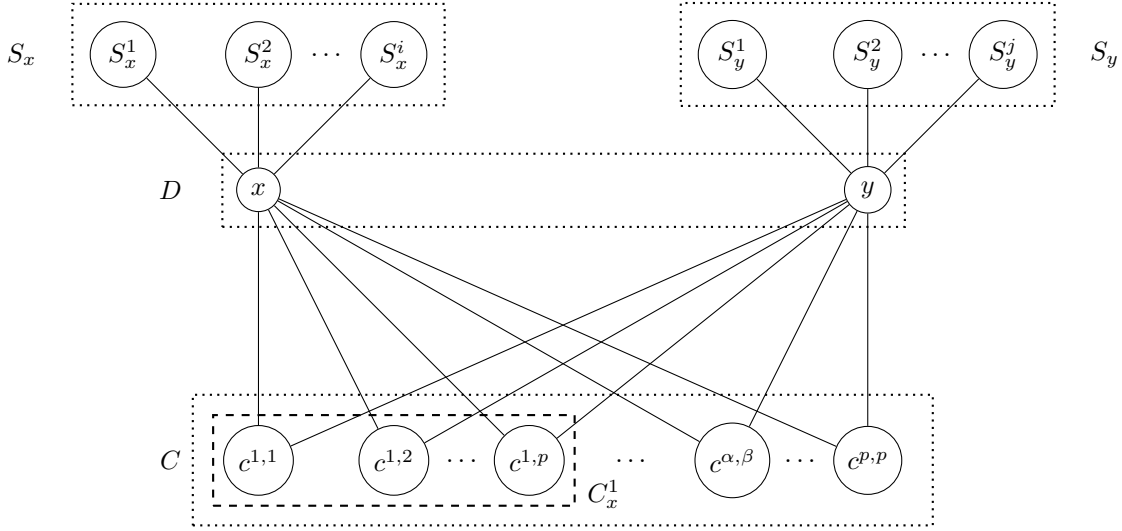
\begin{figure}[t]
		\begin{center}
			\begin{tikzpicture}  
				[scale=.75,auto=center] 
				
				

				\node[circle,draw] (s1) at (-2,8) {$S_x^{1}$};
                \node[circle,draw] (s2) at (0,8) {$S_x^{2}$};
				\node[circle,draw] (s3) at (2,8) {$S_x^{i}$};
				\node[circle,draw]  (s4) at (11,8) {$S_y^{j}$};
				\node[circle,draw]  (s5) at (9,8) {$S_y^{2}$};
                \node[circle,draw] (s6) at (7,8) {$S_y^{1}$};
				\node[circle,draw] (x) at (0,6) {$x$};
                \node[circle,draw] (y) at (9,6) {$y$};
				\node[circle,draw] (c1) at (0,2) {$c^{1,1}$};
				\node[circle,draw]  (c2) at (2,2) {$c^{1,2}$};
				\node[circle,draw]  (c3) at (4,2) {$c^{1,p}$};
                \node[circle,draw] (c4) at (7,2) {$c^{\alpha,\beta}$};
				\node[circle,draw]  (c5) at (9,2) {$c^{p,p}$};
				\node at (3,2) {$\ldots$};
                \node at (5.5,2) {$\ldots$};
                \node at (8,2) {$\ldots$};
                 \node at (1,8) {$\ldots$};
                  \node at (10,8) {$\ldots$};
                  \node at (-3.5,8) {$S_x$};
                    \node at (12.5,8) {$S_y$};
                    \node at (-1.3,6) {$D$};
                    \node at (-1.3,2) {$C$};
                     \node at (5.1,1.5) {$C^{1}_x$};

                \draw (x) -- (c1);
                \draw (x) -- (c2);
                \draw (x) -- (c3);
                \draw (x) -- (c4);
                \draw (x) -- (c5);
                 \draw (y) -- (c1);
                \draw (y) -- (c2);
                \draw (y) -- (c3);
                \draw (y) -- (c4);
                \draw (y) -- (c5);
                 \draw (x) -- (s1);
                 \draw (x) -- (s2);
                \draw (x) -- (s3);
                \draw (y) -- (s4);
                \draw (y) -- (s5);
                \draw (y) -- (s6);
				\draw[thick,dotted] ($(s1.north west)+(-0.4,0.4)$)  rectangle ($(s3.south east)+(0.4,-0.4)$);
            \draw[thick,dotted] ($(s6.north west)+(-0.4,0.4)$)  rectangle ($(s4.south east)+(0.4,-0.4)$);
             \draw[thick,dotted] ($(x.north west)+(-0.3,0.3)$)  rectangle ($(y.south east)+(0.3,-0.3)$);
            \draw[thick,dotted] ($(c1.north west)+(-0.6,0.6)$)  rectangle ($(c5.south east)+(0.6,-0.6)$);
             \draw[thick,dashed] ($(c1.north west)+(-0.3,0.3)$)  rectangle ($(c3.south east)+(0.3,-0.3)$);

			\end{tikzpicture}  
			
		\end{center}
		\caption{A schematic diagram presenting the notation for $H$.}
		\label{fig private neighbors}
	\end{figure}

Let $S = V(H) \setminus [C \cup D]$ denote the set of all private neighbors of $x$ and that of $y$. Moreover, let $S_x = N(x) \setminus C$ and $S^{\alpha}_x = N^{\alpha}(x) \setminus C$ denote the set of private neighbors and private $\alpha$-neighbors of $x$, respectively. Similarly, let 
$S_y = N(y) \setminus C$ and 
$S^{\beta}_y = N^{\alpha}(y) \setminus C$ denote the set of private neighbors and private $\beta$-neighbors of $y$, respectively. See Fig.~\ref{fig private neighbors} for a schematic diagram presenting the proof specific notation.

\begin{remark}\label{rem proof case}
Theorem~1.7 of~\cite{nandy2016outerplanar} proves Theorem~\ref{th main} for $(n,m)=(1,0)$. The proof strategy of Theorem~1.7. from~\cite{nandy2016outerplanar} involves separate treatments depending on the size of $|C|$. Some parts of the proof assuming $|C| \leq 5$  only use the fact that each vertex of $S_x \cup \{x\}$ is connected by an arc/edge or a $2$-path to each vertex of $S_y \cup \{y\}$ and draws contradiction to the maximality of $D$. Those parts of the proof also hold in our setup. Therefore, we can use some observations, that are implicitly proved in~\cite{nandy2016outerplanar}, in our proof. On the other hand, the proofs assuming $|C| \geq 6$ cannot be extended to our set-up, hence the main contribution of this article is in finding a proof under this assumption. 
\end{remark}

\subsubsection*{Building tools for the proof}  

We present some important observations which will be repeatedly used for the proof.

\begin{observation}\label{obs at most 3}
    If two distinct vertices  $u$ and $v$ of $H$ are adjacent or have at least $6$ common neighbors, 
        then $|N^{\alpha}(u) \cap N^{\beta}(v)| \leq 3$. 
\end{observation}

\begin{proof}
    Note that the vertices of $|N^{\alpha}(u) \cap N^{\beta}(v)|$ cannot see each other through $u$ or $v$. 
    Thus, if $|N^{\alpha}(u) \cap N^{\beta}(v)| \geq 4$, then they cannot see each other keeping $H$ planar. 
\end{proof}

\begin{observation}\label{obs adjacent}
For any $\alpha, \beta \in A_{n,m}$ the following are true.
\begin{enumerate}[(i)]
\item If $(u,v) \in S_x \times S_y$ or 
$S^{\alpha}_x \times S^{\alpha}_x$ or 
$S^{\beta}_y \times S^{\beta}_y$, 
then $u$ and $v$ must belong 
to the same or adjacent regions. 

\item If 
$(u,v) \in S^{\alpha}_x \times C^{\alpha}_x$ or 
$S^{\beta}_y \times C^{\beta}_y$ and $u$ belongs to the region $R$, 
then $v$ must belong to the boundary of either $R$ or its adjacent region. 
\end{enumerate}
\end{observation}

\begin{lemma}[~\cite{nandy2016outerplanar}]\label{lem nandy}
    In the planar graph $H$, the following are true when $S_x, S_y \neq \emptyset$. 
    \begin{enumerate}[(i)]
        \item If $|C| \leq 2$, 
then $|V(H)| \leq 3(2n+m)^2+(2n+m)+1$. 

        \item If $|C| =3$ and the vertices of $S$ are spread across all $3$ regions, then $|V(H)| \leq 3(2n+m)^2+(2n+m)+1$. 

        \item If $|C| =4$ and the vertices of $S$ are spread across all $4$ regions, 
        then $|V(H)| \leq 3(2n+m)^2+(2n+m)+1$.       
    \end{enumerate}
\end{lemma}

\begin{proof}
    The proofs are implicitly presented in~\cite{nandy2016outerplanar} within the proofs of Lemmas~4.6, 4.7, and 4.8 (Claims~3 and~4).  
\end{proof}

A direct consequence of Observation~\ref{obs adjacent} is the following. 

\begin{observation}\label{obs spread region}
In the planar graph $H$, the following are true. 
    \begin{enumerate}[(i)]
        \item If $|C| \geq 4$, then the vertices of 
        $S^{\alpha}_x \cup S^{\beta}_y$ can be spread across at most $2$ regions, and those $2$ regions must be adjacent.

 \item If $|C| \geq 3$, and $S_x, S_y \neq \emptyset$, then the vertices of $S$ can be spread across at most three regions of the type $R_{i-1}, R_i$ and $R_{i+1}$, except for the case described in 
 Lemma~\ref{lemma C-S trade off}(iii).   
    \end{enumerate}
\end{observation}

The next lemma discusses the relation between the size of $C$ and that of $S$.

\begin{lemma}\label{lemma C-S trade off}
In the planar graph $H$, the following are true. 
    \begin{enumerate}[(i)]
\item If $|C^{\alpha}_t| \geq 5$, then 
        $|S^{\alpha}_t| =0$, where $t \in \{x,y\}$.

\item If $|C^{\alpha}_t| \geq 4$, then 
        $|S^{\alpha}_t| \leq 2$, where $t \in \{x,y\}$. Moreover, if $|C| \geq 5$, then 
$|S^{\alpha}_t| \leq 1$.

\item If $|C| \geq 3$ and the vertices of $S$ are not spread across all the regions, then 
$|S^{\alpha}_x \cup S^{\beta}_y| \leq 3(2n+m)+1$.

\item If $|C| \geq 3$ and 
$S_x, S_y \neq \emptyset$, then 
$$|S^{\alpha}_x| \leq 3(2n+m)+1 - 
              \max_{\gamma}\{|S^{\gamma}_y|\}
\text{ and } 
|S^{\beta}_y| \leq 3(2n+m)+1 - 
              \max_{\gamma}\{|S^{\gamma}_x|\}.$$ 
In particular, 
$|S^{\alpha}_x|, |S^{\beta}_y| \leq 3(2n+m)$. 
\end{enumerate}
\end{lemma}

\begin{proof}
    Let $x_1$ be a vertex of $S_x$ in the region $R_2$. The only vertices of $C$ that $x_1$ can see without seeing through $x$ are $c_0$ (through $c_1$), $c_1$, $c_2$, and $c_3$ (through $c_2$).

    (i) Due to the above, one can observe that, a vertex of  $S^{\alpha}_x$ can see at most $4$ vertices of $C^{\alpha}_x$. 

    (ii) Suppose that a vertex $x_1$ of $S^{\alpha}_x$ sees exactly $4$ vertices of $C^{\alpha}_x$. Without loss of generality, if 
    $x_1$ belongs to $R_2$, 
    then the $4$ vertices of $C^{\alpha}_x$ that $x_1$ sees are exactly $c_0, c_1, c_2, c_3$ and $x_1$ is adjacent to $c_1, c_2$. Thus there can be at most one vertex of $S^{\alpha}_x$ in the region $R_2$. If $|C|=4$, then there can be at most one more vertex of $S^{\alpha}_x$ in an adjacent region. If $|C| \geq 5$, then there cannot be any other 
    vertex of $S^{\alpha}_x$ in any other region. 

    (iii)  Without loss of generality assume that there are no vertices of 
    $S^{\alpha}_x \cup S^{\beta}_y$ in the region $R_0$. Now we  modify the graph $H$ by deleting all arcs and edges contained in the region $R_0$, and then adding one edge between the vertices $x,y$ that goes through the region $R_0$. Next, we contract the edge $xy$ (that goes through the region $R_0$) and the new vertex obtained is called $z$. Notice that, this so-obtained graph, say $H'$, is planar and $z$ is a 
    universal vertex (that is, adjacent to every other vertex). 
    Therefore, $H' - z$ is outerplanar. 

    Note that, in $H$, the vertices of $S^{\alpha}_x \cup S^{\beta}_y$ were seeing each other without 
    any contribution from any vertex, arc, or edge contained in region $R_0$. Therefore, by Theorem $3.1$ of~\cite{bensmail2017analogues} , we have $$|S^{\alpha}_x \cup S^{\beta}_y| \leq 3(2n+m)+1.$$ 

\medskip

(iv)  Follows directly from Lemma~\ref{lemma C-S trade off}(iii). 
\end{proof}

\subsubsection*{When $S_x = \emptyset$ or $S_y = \emptyset$}
In this part of the proof, we prove Theorem~\ref{th main} if $x$ or $y$ has no private neighbors.

\begin{lemma}\label{lem 2empty}
    If $S_x$ and $S_y$ both are empty, 
    then $H$ has at most $3(2n+m)^2+(2n+m)+1$ vertices.
\end{lemma}

\begin{proof}
   If $|C| \geq 6$, then we can apply Observation~\ref{obs at most 3} to imply
   $$|C| \leq |\cup_{(\alpha, \beta)} C^{\alpha,\beta}| \leq \sum_{(\alpha, \beta)} |C^{\alpha,\beta}| \leq 3(2n+m)^2.$$
   Thus, as $3(2n+m)^2$ is always greater than $5$, we have 
   $$|V(H)| \leq |D|+|C|+|S| \leq 2+ 
   \max\{5, 3(2n+m)^2\} +0 \leq 3(2n+m)^2+(2n+m)+1,$$
   and we are done. 
   \end{proof}

\begin{lemma}\label{lem 1empty csmall}
    If $|C| \leq 5$, and if exactly one between $S_x$ and $S_y$ is empty, 
    then $H$ has at most $3(2n+m)^2+(2n+m)+1$ vertices.
\end{lemma}

\begin{proof}
    Suppose, without loss of generality, that $S_x \neq \emptyset$ and $S_y = \emptyset$. Note that, $x$ must be non-adjacent to $y$ as, otherwise, $x$ will be a dominating vertex. Notice that the vertices of $S = S_x$ must all be adjacent to some 
    $c_i \in C$ to see $y$.  That means, if $|C| \leq 5$ and  $S$ has at least $16$ vertices, 
    then some $c_i$ will  be adjacent to at least 
    $4$ vertices from $S_x$ by the pigeonhole principle, and $2$ from $C$ due to triangulation. 
    Thus, notice that $\{x,c_i\}$ is a dominating set 
    of $H$ having at least $6$ common neighbors, which contradicts the maximality of $D$. 
    Therefore, if $|C| \leq 5$, then $|S| \leq 15$ and we are done. 
\end{proof}

\begin{lemma}\label{lem 1empty helper}
    If $|C| \geq 6$ and exactly one between $S_x$ 
    and $S_y$ is empty, 
    then $|C^{\alpha}_x| + |S^{\alpha}_x| \leq 3(2n+m)$
    for all $\alpha \in A_{n,m}$. 
\end{lemma}

\begin{proof}
We know that $|C^{\alpha}_x| \leq 3(2n+m)$ due to Observation~\ref{obs at most 3}. Thus, 
by Lemma~\ref{lemma C-S trade off}(i) and~(ii)  we can conclude that, 
if $|C^{\alpha}_x| \geq 4$, then 
    $$|C^{\alpha}_x| + |S^{\alpha}_x| \leq \max\{3(2n+m), 6\} \leq 3(2n+m).$$
    For $|C^{\alpha}| \leq 3$, first we will handle an exceptional case. 
    If $|C^{\alpha}|\leq 2$ and 
    $|S^{\alpha}| \leq 4$, then 
    $$|C^{\alpha}_x| + |S^{\alpha}_x| \leq 6 \leq 3(2n+m).$$
    In the other cases, we first show that all vertices of $S^{\alpha}_x$ are adjacent to a particular vertex of $C$ (say, $c_1$). 
    
    If the vertices of $S^{\alpha}_x$  are spread across two adjacent regions, say $R_1$ and $R_2$, then 
    all vertices of $S^{\alpha}_x$ are adjacent to $c_1$ to see each other. If all vertices of $S^{\alpha}_x$ are contained in only one region, say $R_2$, and $|C^{\alpha}_x| = 3$, then we may assume without loss of generality that $c_0 \in C^{\alpha}_x$, and hence, all vertices of $S^{\alpha}_x$ are adjacent to $c_1$ to see $c_0$. If $|C^{\alpha}_x| \leq 2$, then we have nothing to prove unless $|S^{\alpha}_x| \geq 5$. Thus, let us assume that $|S^{\alpha}_x| \geq 5$.
    As $x$ and $y$ are non-adjacent to avoid $H$ being dominated by $x$, the vertices of $S^{\alpha}_x$ must be adjacent either $c_1$ or $c_2$ to see $y$. Thus, without loss of generality assume that $3$ vertices of $S^{\alpha}_x$ are adjacent to $c_1$. This will force all vertices of $S^{\alpha}_x$ to be adjacent to $c_1$ for seeing the $1$-nearest vertex of $S^{\alpha}_x$ to $xc_1$.

Assume that $c_0, c_2, y$ are $\beta, \gamma, \tau$-neighbors of $c_1$, respectively. 
Let $Z_1$ and $Z_2$ be the set of vertices of $S^{\alpha}_x$ that are contained in the regions $R_1$ and $R_2$, respectively. Furthermore, let 
$z_{i1}$ and $z_{i2}$ denote the $1$-farthest and the $2$-farthest vertices of $Z_i$ to $xc_1$, respectively, for all $i \in \{1,2\}$.  
Notice that, the vertices from 
$Z_1 \cup Z_2 \setminus \{z_{11}, z_{21}\}$ sees $y$ through $c_1$, and thus none of them are $\tau$-neighbors of $c_1$. Moreover, if $c_0 \in C^{\alpha}_x$, 
then the vertices of 
$Z_1 \cup Z_2 \setminus \{z_{11}, z_{12}\}$   
sees $c_0$ 
through $c_1$, and thus none of them are $\beta$-neighbors  of $c_1$.
Similarly, if  $c_2 \in C^{\alpha}_x$, 
then the vertices of 
 $Z_1 \cup Z_2 \setminus \{z_{21}, z_{22}\}$
sees  $c_2$ through $c_1$, and thus none of them are  $\gamma$-neighbors of $c_1$.

\medskip

$\bullet$ If $|C^{\alpha}_x|=3$ and $\beta = \gamma$, then $S^{\alpha}_x$ does not contain any $\beta$-neighbor of $c_1$. 
Therefore, by Observation~\ref{obs at most 3} we have
$$|C^{\alpha}_x| + |S^{\alpha}_x| \leq 3 + 3(2n+m-1) = 3(2n+m).$$

$\bullet$ If $|C^{\alpha}_x|=3$ and $\beta \neq \gamma$, 
then $S^{\alpha}_x \setminus \{z_{11}, z_{12}, z_{21}\}$ does not contain any $\beta$-neighbor or any $\tau$-neighbor of $c_2$. 
Therefore, by Observation~\ref{obs at most 3} we have
$$|C^{\alpha}_x| + |S^{\alpha}_x \setminus \{z_{11}, z_{12}, z_{21}\}| + |\{z_{11}, z_{12}, z_{21}\}| \leq 3 + 3(2n+m-2) + 3 = 3(2n+m).$$

\medskip

  $\bullet$  If $|C^{\alpha}_x|=2$ and $\beta = \tau$, then without loss of generality assume that $C^{\alpha}_x = \{c_0, c_1\}$. 
  In this case,  
 $S^{\alpha}_x \setminus \{z_{11}\}$ does not contain any $\beta$-neighbor of $c_1$. Therefore, by Observation~\ref{obs at most 3} we have
$$|C^{\alpha}_x| + |S^{\alpha}_x \setminus \{z_{11}\}| +|\{z_{11}\}| \leq 2 + 3(2n+m-1) +1 = 3(2n+m).$$

$\bullet$ If $|C^{\alpha}_x|=2$ and $\beta \neq \tau$,  then assuming $C^{\alpha}_x = \{c_0, c_1\}$ we can note that $S^{\alpha}_x \setminus \{z_{11}, z_{12}\}$ does not contain any $\beta$-neighbor or any $\tau$-neighbor of $c_1$. 
Therefore, by Observation~\ref{obs at most 3} we have
$$|C^{\alpha}_x| + |S^{\alpha}_x \setminus \{z_{11}, z_{12}\}| + |\{z_{11}, z_{12}\}| \leq 2 + 3(2n+m-2) + 2 \leq 3(2n+m).$$

\medskip

$\bullet$  If $|C^{\alpha}_x|=1$,  
then $S^{\alpha}_x \setminus \{z_{11}, z_{21}\}$ does not contain any $\tau$-neighbor of $c_1$. Therefore, by Observation~\ref{obs at most 3} we have
$$|C^{\alpha}_x| + |S^{\alpha}_x \setminus \{z_{11}, z_{21}\}| + |\{z_{11}, z_{21}\}| \leq 1 + 3(2n+m-1) + 2 = 3(2n+m).$$

\medskip

This completes the proof. 
\end{proof}

The following lemma is a crucial observation using Lemmas~\ref{lem 2empty},\ref{lem 1empty csmall}, and~\ref{lem 1empty helper}.

\begin{lemma}\label{lem 1empty final}
    If $|C| \geq 6$ and at least one between $S_x$ 
    and $S_y$ is empty, 
    then $H$ has at most $3(2n+m)^2+(2n+m)+1$ vertices. 
\end{lemma}

\begin{proof}
    We know that the statement is true if $S_x = S_y = \emptyset$ and if $|C| \leq 5$ due to Lemmas~\ref{lem 2empty} and~\ref{lem 1empty csmall}, respectively.

    If $|C| \geq 6$, and if exactly one between $S_x$ and $S_y$ is non-empty, then due to Lemma~\ref{lem 1empty helper} we have $|C^{\alpha}_x| + |S^{\alpha}_x| \leq 3(2n+m)$ for all 
    $\alpha \in A_{n,m}$. 
    Therefore, 
    $$|V(H)| \leq |D| + |C| + |S| \leq 2 + \sum_{\alpha} (|C^{\alpha}_x| + |S^{\alpha}_x|) \leq 2 + 3(2n+m)^2 \leq 3(2n+m)^2 + (2n+m) + 1.$$
    This completes the proof of the lemma. 
\end{proof}

\subsubsection*{When $S_x \neq \emptyset$ and $S_y \neq \emptyset$}
Let us assume that $x$ has $i$ different types 
and $y$ has $j$ different types of neighbors in $S$. 
Without loss of generality, let us suppose that 
$i \geq j$ for the rest of the section. 
A common neighbor $c \in C^{\alpha, \beta}$ of $x,y$ is an
\textit{excess} if $S^{\alpha}_x \cup S^{\beta}_y \neq \emptyset$. The set of all excesses is denoted by $E$. 
The lemmas in this phase of the proof will assume the above-mentioned setup.

\begin{lemma}\label{lem key bigC}
If $|C| \geq 3$, $S_x$ and $S_y$ are both 
non-empty, and $i \geq j$,
    then we have
    $$|V(H)| \leq [3(2n+m)^2 + (2n+m) + 1] 
    - [3j(2n+m - i) + (2n+m) +(i-j)s_{max} -j -|E| -1]$$
    where $s_{max} = \displaystyle{\max_{\gamma}}\{|S^{\gamma}_y|\}$.
\end{lemma}

\begin{proof}
    Due to Observation~\ref{obs spread region}, 
    without loss of generality,
    we may assume that the vertices of 
    $S$ are contained in the 
    regions $R_1, R_2, R_3$, and that, $R_2$ contains some vertices of $S$. 
    Note that, if a vertex $u \in S$ sees a vertex $v\in C$ without the help of $x$ or $y$, then 
    we are forced to have $v \in \{c_{k-1}, c_0, c_1, c_2, c_3, c_4\}$. Notice that, each vertex of the excess set $E$ must be seen by some vertex of $S$ without the help of $x$ or $y$. Thus, 
    $E \subseteq \{c_{k-1}, c_0, c_1, c_2, c_3, c_4\}$, which implies $|E| \leq 6$. However, this allows us to use Observation~\ref{obs at most 3} to estimate the number of vertices in 
    $C \setminus E$. 

    We know that the vertices of $C \setminus E$ have at most $(2n+m-i)$ types of neighbors of $x$ and at most $(2n+m-j)$ different types of neighbors of $y$. Thus, by Observation~\ref{obs at most 3} 
    \begin{equation}\label{equation estimate C}
    |C \setminus E| \leq 3(2n+m-i)(2n+m-j) = 3(2n+m)^2 - 3(2n+m)(i+j) +3ij.    
    \end{equation}
    On the other hand, if we assume $i \geq j$ without loss of generality, then it is possible to find $j$ distinct sets of the type $S^{\alpha}_x \cup S^{\beta}_y$ and $(i-j)$ sets of the type $S^{\gamma}_x$ which partitions $S$. 
Using the bounds from Lemma~\ref{lemma C-S trade off}(iii) and (iv), we have 
\begin{equation}\label{equation estimate S}
    |S| \leq [3(2n+m)+1]j + [3(2n+m) - s_{max}](i-j) = 3(2n+m)i + j - (i-j)s_{max}
\end{equation}
where $s_{max} = \displaystyle{\max_{\gamma}}\{|S^{\gamma}_y|\}$. 
    Thus, combining Equations~(\ref{equation estimate C}) and~(\ref{equation estimate S}) we get the following bound:
\begin{align*}
    |V(H)| &\leq |D| + |C| + |S| \\
    &\leq 2 + 3(2n+m)^2 - 3(2n+m)(i+j) +3ij + |E| + 3(2n+m)i + j - (i-j)s_{max}\\
    &= 3(2n+m)^2 - 3j(2n+m - i) + j + |E| + 2 - (i-j)s_{max}\\
    &= [3(2n+m)^2 + (2n+m) + 1] 
    - [3j(2n+m - i) + (2n+m) + (i-j)s_{max} -j -|E| -1].   
\end{align*}
This completes the proof of the lemma. 
\end{proof}

A direct corollary using the equation established is the following. 

\begin{corollary}\label{cor easy case}
Let $|C| \geq 3$, $S_x$ and $S_y$ are both 
non-empty, and $i \geq j$. Moreover, let $H$ satisfy one of the following three conditions: 
\begin{itemize}
    \item $i \leq (2n+m)-2$,

    \item $i = (2n+m) - 1$ and $j \geq 2$, 

    \item $i = (2n+m) - 1$, $j = 1$, and $|E| \leq 4$.
\end{itemize} 
   Then  $H$ has at most $3(2n+m)^2+(2n+m)+1$ vertices. 
\end{corollary}

\begin{proof}
    Substituting each of the listed values from the statement of this corollary in Lemma~\ref{lem key bigC} implies the proof. 
\end{proof}

Using Lemma~\ref{lem key bigC} we will prove the theorem for the cases of $|C| \geq 3$ which did not get captured by Corollary~\ref{cor easy case}.

\begin{lemma}\label{lem p-1,1}
   Let $|C| \geq 3$, $S_x$ and $S_y$ are both 
non-empty,   $i = (2n+m)-1$, $j=1$, and 
$|E| \geq 5$. Then  $H$ has at most $3(2n+m)^2+(2n+m)+1$ vertices. 
\end{lemma}

\begin{proof}
First note that, if 
$\displaystyle{\max_{\gamma}}\{|S^{\gamma}_y|\} \geq 2$, 
then by Lemma~\ref{lem key bigC}, as $2n+m \geq 3$, we have 
$$|V(H)| \leq [3(2n+m)^2 + (2n+m) + 1] 
    - [3 + 3 + 2 -1 - 6  - 1] = 
    3(2n+m)^2 + (2n+m) + 1.$$
Thus, we may assume that 
$S = S^{\beta}_y = \{y_1\}$. 
In this case, Lemma~\ref{lem key bigC}   gives us 
$$|V(H)| \leq [3(2n+m)^2 + (2n+m) + 1] 
    - [3 + 3 + 1 -1 - 6  - 1] = 
    3(2n+m)^2 + (2n+m) + 2.$$
    That means, unless $|S^{\alpha}_x| = 3(2n+m)$
    for all $\alpha$ where $S^{\alpha}_x$ is non-empty and $|E|=6$, we are done. 
    Thus, for the next parts of the proof, we can use the additional information that 
    if $S^{\alpha}_x \neq \emptyset$, 
    then $|S^{\alpha}_x| = 3(2n+m) \geq 9$.

As $|E| = 5$, we have $c_{k-1} \in E$, and we may assume that  
$c_{k-1} \in C^{\alpha_1, \beta_1}$ 
for some $\alpha_1, \beta_1$.  
Thus, either 
$S^{\alpha_1}_x \neq \emptyset$ or 
$S^{\beta_1}_y \neq \emptyset$. 
If $S^{\alpha_1}_x \neq \emptyset$, then all its vertices are contained in $R_1$ and they must see $c_{k-1}$  through $c_0$. In such a scenario, $y_1$ must see $x_1$ through $c_0$, where $x_1$ is the vertex of $S^{\alpha_1}_x$ which is $1$-nearest to $xc_0$. That is, $y_1$ must be contained in $R_1$. 
If $S^{\beta_1}_y \neq \emptyset$, 
then the only vertex in it, that is, $y_1$ must belong to $R_1$ and see $c_{k-1}$ through $c_0$. Thus, we have established that  $y_1$ is in $R_1$, which implies that $R_3$ does not contain any vertex of $S$, due to Observation~\ref{obs spread region} and the fact that 
$S_y = \{y_1\}$. Hence, $c_4$ cannot be a vertex of $E$, and thus, $|E| \leq 5$, a contradiction. 
\end{proof}

\begin{lemma}\label{lem p notp}
   Let $|C| \geq 3$, $S_x$ and $S_y$ are both 
non-empty, and  $j < i = (2n+m)$. Then  $H$ has at most $3(2n+m)^2+(2n+m)+1$ vertices. 
\end{lemma}

\begin{proof}
As $i = 2n+m$, we must have $C=E$. 
Thus, $|C| = |E| \leq 6$. 
Suppose $|S^{\alpha}_x| \leq 4$ for some $\gamma$. 
Notice that, we have obtained the inequality in Lemma~\ref{lem key bigC} by estimating $|S^{\alpha}_x| = 3(2n+m)$. 
Thus, with the additional information of 
$|S^{\alpha}_x| \leq 4$, we can obtain the following 
modified inequality:
\begin{align*}
|V(H)| &\leq [3(2n+m)^2 + (2n+m) +1] - [3(2n+m)+(i-j)s_{max} - j -|E| -1] - [(2n+m)-4]\\
&\leq [3(2n+m)^2 + (2n+m) +1] - [2 -|E| -1] - 5\\
&\leq [3(2n+m)^2 + (2n+m) +1].
\end{align*}
Thus, we may assume that $|S^{\alpha}_x| \geq 5$ for all $\alpha$. 
Therefore, it is not possible for the vertices of 
$S^{\alpha}_x$ to be spread across $3$ different regions due to Observation~\ref{obs spread region}. 

Suppose without loss of generality, one vertex $x_1$ of $S^{\alpha}_x$ belongs to a region, say, $R_2$, and there exists a $\gamma$-neighbor of $x$ from $C$, say $c_0$, which is not on the boundary of $R_2$. 
Then without loss of generality, all the vertices of 
$S^{\alpha}_x$ must be adjacent to $c_1$ to see $c_0$ or each other. However, this will force all the vertices of $S_y$ to be adjacent to $c_1$ to see the vertices of $S^{\alpha}_x$. Thus, $D' = \{x, c_2\}$ 
is a dominating set with 
$S^{\alpha}_x \cup \{c_0\} \subseteq N(x) \cap N(c_1)$ implying $|N(x) \cap N(c_1)| \geq 6$. Hence, unless $|E| = 6$, the existence of $D'$ contradicts the maximality of $D$. 

Thus, now we may assume that $|E| = 6$. In particular, we may assume that 
$C= E = \{c_{k-1}$, $c_0$, $c_1$, $c_2$, $c_3$, $c_4\}$ and that 
the vertices of $S$ are contained in the $3$ regions $R_1, R_2, R_3$. Suppose $c_{k-1}$, and $c_4$ are a 
$\alpha_1$-neighbor and $\alpha_2$-neighbor of $x$, respectively. 
This will imply that all the vertices of 
$S^{\alpha_1}_x$ are contained in $R_1$ and 
sees $c_{k-1}$ through $c_0$.
 As $S^{\alpha_1}_x$ has at least $5$ vertices, this will force all the vertices of $S_y$ to see the $1$-nearest vertex of $S^{\alpha_1}_x$ to $xc_0$ through $c_0$. In particular, all the vertices of 
 $S_y$ must belong to $R_1$. Therefore, there cannot be any vertex of $S$ in $R_3$, which contradicts the fact that $c_4 \in E$, by the definition of $E$. 
 Hence we are done. 
\end{proof}

\begin{lemma}\label{lem p p}
   Let $|C| \geq 3$, $S_x$ and $S_y$ are both 
non-empty, and  $ i = j = (2n+m)$. Then  $H$ has at most $3(2n+m)^2+(2n+m)+1$ vertices. 
\end{lemma}

\begin{proof}
If $S^{\alpha}_x$ are spread across exactly $2$ 
adjacent regions, then $|E| \leq 5$ and  
all vertices of $S^{\alpha}_x$ are adjacent to some vertex of $C$, say $c_1$, forcing all vertices of $S_y$ to be adjacent to $c_1$. As $S_y$ has a least $3$ vertices (since $j = 2n+m \geq 3$), all vertices of $S_x$ are also adjacent to $c_1$. Thus, both 
$D' = \{x, c_1\}$ and $D'' = \{y, c_1\}$ are dominating set. Due to the maximality of $D$, we can infer that 
    \begin{equation}
        |S_x| \leq |N(x) \cap N(c_1)| \leq |C| = |E| \leq 5 \text{ and } |S_y| \leq |N(y) \cap N(c_1)| \leq |C| = |E| \leq 5.
    \end{equation}
    Therefore, 
    \begin{equation}
        |V(H)| \leq |D|+|C|+|S| \leq 2+|E| + |S_x| +|S_y| \leq 2+5+5+5 =17.
    \end{equation}
    This means, if all vertices of $S$ are adjacent to a single vertex $c_1$ (say) of $C$, then we are done. We will use this information multiple times in the proof of this lemma. On the other hand, for the rest of the proof, we may now assume that  $S^{\alpha}_x$ (resp., $S^{\beta}_y$) 
    is part of a single region for all $\alpha \in A_{n,m}$ (resp., $\beta \in A_{n,m}$).

 Note that the inequality due to Lemma~\ref{lem key bigC} modifies to
    \begin{equation}\label{equation i=j=p}
        |V(H)| \leq [3(2n+m)^2+(2n+m)+1] + |E| +1
    \end{equation}
    where every set of the form 
    $S^{\alpha}_x \cup S^{\beta}_y$ is estimated to have exactly $3(2n+m)+1 \geq 10$ vertices. Thus, 
    in particular, we are done unless 
    \begin{equation}\label{eqn s size}
        |S^{\alpha}_x \cup S^{\beta}_y| \geq 10-|E| 
    \end{equation}
    for all $\alpha, \beta \in A_{n,m}$.

If $S$ is spread across $2$ adjacent regions, then $|E| \leq 5$ and we must have a situation where the vertices of $S^{\alpha}_x$ and $S^{\beta}_y$ are contained in $2$ adjacent regions, say $R_1$ and $R_2$ without loss of generality. Thus, the vertices of $S^{\alpha}_x \cup S^{\beta}_y$ must all be adjacent to $c_1$ to see each other. As 
$$|S^{\alpha}_x \cup S^{\beta}_y| \geq 10 -|E| \geq 10-5 = 5,$$
we must have at least $3$ vertices in $S^{\alpha}_x$ without loss of generality. This forces all vertices of $S_y$ to be adjacent to $c_1$, which in turn implies all vertices of $S_x$ to be adjacent to $c_1$, and thus we are done.

    If all vertices of 
    $S$ belong to a single region, 
    say $R_2$, and $c_0 \in C^{\alpha, \beta}$, 
    then  $|E| \leq 4$ and the vertices of 
    $S^{\alpha}_x \cup S^{\beta}_y$ must see $c_0$ through $c_1$, without loss of generality.  
    As due to inequality~(\ref{eqn s size}) 
    $$|S^{\alpha}_x \cup S^{\beta}_y| \geq 10-|E| \geq 10-4=6.$$ 
    Thus, without loss of generality, we may assume that
    at least $3$ vertices of $S^{\alpha}_x$  are adjacent to $c_1$. This will force all vertices of $S_y$ to be adjacent to $c_1$, which in turn implies all vertices of $S_x$ to be adjacent to $c_1$, and thus we are done.   
\end{proof}
\textbf{Concluding the proof: } The proof of Theorem~\ref{th main} directly follows from Lemma~\ref{lem 1empty final}, Corollary~\ref{cor easy case}, and Lemmas~\ref{lem p-1,1}, \ref{lem p notp}, \ref{lem p p}.

	\bibliographystyle{abbrv}
	\bibliography{mybibliography}

\end{document}